# REGIONAL BOUNDARY EXPONENTIAL REDUCED OBSERVABILITIY IN DISTRIBUTED PARAMETER SYSTEMS

## Raheam A. Al-Saphory[1*] and Hind K. Kolaib[1]

[1]Department of Mathematics, College of Education for Pure Sciences, Tikrit University, Tikrit, Iraq.

[1*]E-mail: saphory@hotmail.com and [1]Hind.K.Kolaib@hotmail.com

*To whom correspondence should be addressed*

**ABSRACT:** The aim of this chapter is to introduce the concept of regional boundary exponential reduced observability in connection with the sensors characterizations on a sub-region $\Gamma$ of the considered domain boundary $\partial\Omega$. More precisely, we explore the original results devoted to this concept in linear dynamical systems which is generated by a strongly continuous semi-group in Hilbert space $H^1(\Omega)$. Thus, the existence of sufficient conditions is presented and examined for regional boundary exponential reduced estimator in parabolic infinite dimensional systems. Finally, we apply these results to the exchange systems with various strategic sensors.

**Keywords:** $\Gamma_{E\mathcal{R}}$-observability, $\Gamma_E$-detectability, $\Gamma_{E\mathcal{R}}$-strategic sensor, $\Gamma_{E\mathcal{R}}$- detectability, exchange system, $\Gamma_{E\mathcal{R}}$-strategic sensor.

**2010 AMS SUBJECT CLASSIFICATION: 93A30; 93B07; 93C05; 93C20**

## 1. INTRODUCTION

For a distributed parameter system evolving on a special domain $\Omega$, the observability concept has been widely developed and survey of these developments can be found [1-3]. The purpose of an exponential estimator is to provide an exponential state estimation for the considered system state [4]. New direction of regional analysis for infinite dimensional systems has been recently explored by Al-Saphory and El Jai *et al.* in infinite time(may be called regional asymptotic or exponential analysis) as in ref.s [5-11] and for finite time (may be called regional observability) as in [12-14]. In this paper, we introduce and study the notion of exponential regional boundary reduced observability in a given region $\Gamma$ of the domain boundary $\partial\Omega$. Thus the developed approach is an extension of previous works to the regional case as in [6-15]. Moreover the relationship between this notion, regional boundary detectability and regional boundary strategic sensors are tackled and analyzed. The reason behind the study of this notion, there exist some problem in the real world cannot observe the system state in the whole domain boundary, but in a part of this domain boundary. The scenario described by energy exchange problem, where the aim is to determine the energy exchanged in a casting plasma on a plane target which is perpendicular to the direction of the flow from measurements (internal pointwise sensors) carried out by thermocouples [16-17] (Figure 1), where (1) is the torch of plasma, (2) is the probe of (steal), (3) is the insulator, $\Gamma$ is the face of exchange and $b_1$, $b_2$ sensor locations. This paper is organized as follows. Section 2 is devoted to the introduction of regional exponential detectability and considered system with $\Gamma_E$-detectability and $\Gamma$-observability. We study the links of this notion with the regional exponential observability and strategic sensors. In Section 3, we study a regional exponential observability through the relations between $\Gamma_E$-estimator reconstruction method and $\Gamma_E$- observability. In section 4 an section 5 we introduce regional exponential reduced observability notion for a distributed parameter system in terms of regional exponential reduced detectability and reduced strategic sensors. In the last section, we illustrate applications with different domains and circular strategic sensors of two-phase exchange systems.

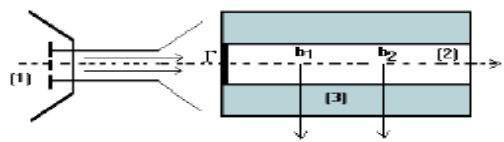

**Fig. 1:** Model of energy exchanged problem on $\Gamma$

## 2. Problem Formulation

Let $\Omega$ be a regular, bounded and open subset of $\mathbb{R}^n$, with boundary $\partial\Omega$. Let $\Gamma$ be a non-empty given sub-region of $\partial\Omega$, with positive measurement. We denote $Q = \Omega \times (0,\infty)$ and $\Theta = \partial\Omega \times (0,\infty)$. Let $X, U$, and $\mathcal{O}$ be a separable Hilbert spaces, where $X$ is the state space, $U$ is the control space and $\mathcal{O}$ the observation space. We consider $X = H^1(\Omega)$, $U = L^2(0,\infty,\mathbb{R}^p)$ and $\mathcal{O} = L^2(0,\infty,\mathbb{R}^q)$ where $p$ and $q$ hold for the numbers of actuators and sensors [18]. We consider the system described by the following parabolic partial deferential equations:

$$\begin{cases} \frac{\partial x}{\partial t}(\xi,t) = Ax(\xi,t) + Bu(t) & Q \\ x(\xi,0) = x_0(\xi) & \overline{\Omega} \\ \frac{\partial x}{\partial v_A}(\eta,t) = 0 & \Sigma \end{cases} \qquad (1)$$

augmented with the output function

$$y(.,t) = Cx(.,t) \qquad (2)$$

Where $\overline{\Omega}$ holds for the closure of $\Omega$ and $x_0(\xi)$ is supposed to be unknown in the state space $X = H^1(\overline{\Omega})$. The system (1) is defined with a Neumann boundary conditions, $\partial x/\partial v_A$ holds for the outward normal derivative. Thus, $A$ is a second-order linear differential operator, and is self-adjoint with compact resolvent. The operators $B \in L(\mathbb{R}^p, X)$ and $C \in L(H^1(\overline{\Omega}), \mathbb{R}^q)$ depend on the strutures of actuators and sensors [19-20]. That means, in the case of pointwise (internal or boundary) and boundary zone sensors (actuators), we have $B \notin L(\mathbb{R}^p, X)$ and $C \notin L(H^1(\Omega), \mathbb{R}^q)$, we refer to see [21]. Thus, the system (1) has a unique solution [22] given by

$$x(\xi,t) = S_A(t)x_0(\xi) + \int_0^t S_A(t-\tau)Bu(\tau)d\tau. \qquad (3)$$

The problem is that how to give an approach which enables to estimate the system state in a sub-region $\Gamma$ of the boundary as in (Figure 2), using convenient sensors. Then, the mathematical model in (Figure 2) is more general spatial case in (Figure 1).

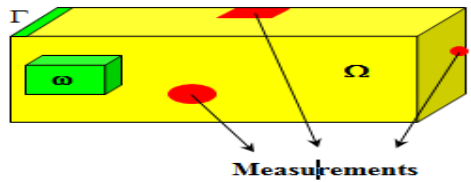

**Fig. 2:** The domain of $\Omega$, the sub-regions $\omega$ and $\Gamma$, various sensors locations.

The regional boundary exponential reduced estimator is defined when the output gives a part of the state vector in this region. We need to rewrite the hypothesizes and some definitions in the following forms:

▪ The operator $K$ is defined by following

$$\begin{aligned} K: X &\to \mathcal{O} \\ x &\to CS_A(.)x \end{aligned}$$

then, we obtain

$$y(.t) = K(t)x(.,0)$$

where $K$ is bounded linear operator (this is valuable on some output function) [21].

▪ The operator $K^*: \mathcal{O} \to X$ is the adjoint of $K$ defined by

$$K^*y^* = \int_0^t S_A^*(s) C^* y^*(.,s)ds$$

- The trace operator of order zero

$$\gamma_0 : H^1(\Omega) \to H^{1/2}(\partial\Omega)$$

is linear, subjective, and continuous [23], such that $x_0^\Gamma$ is the restriction of the trace of the initial state $x_0$ to $\Gamma$. $\gamma_0^*$ denote the adjoint of the operator $\gamma_0$ given by

$$\gamma_0^* : H^{1/2}(\partial\Omega) \to H^1(\Omega)$$

- For a sub-region $\Gamma \subset \partial\Omega$ let $\chi_\Gamma$ be the restriction function defined by

$$\chi_\Gamma : H^{1/2}(\partial\Omega) \to H^{1/2}(\Gamma)$$
$$x \to \chi_\Gamma x = x|_\Gamma$$

where $x|_\Gamma$ is the restriction of $x$ to $\Gamma$. We denote by $\chi_\Gamma^*$ the adjoint function of $\chi_\Gamma$ and defined by

$$\chi_\Gamma^* : H^{1/2}(\Gamma) \to H^{1/2}(\partial\Omega)$$

- Let $\chi_\omega$ be the function defined by

$$\chi_\omega : H^1(\Omega) \to H^1(\omega)$$
$$x \to \chi_\omega x = x|_\omega$$

where $x|_\omega$ is the restriction of the state $x$ to $\omega$.

- The system (1) is $\omega_E$-stable if the operator $A$ generates a strongly continuous semi-group $(S_A(t))_{t \geq 0}$ which is $\omega_E$-stable.

- If The system (1) is $\omega_E$-stable, then the solution of autonomous system associated with (1)-(2), converges exponentially to zero when $t$ tends to $\infty$

- The systems (1)-(2) are said to be regionally exponential detectable on $\omega$ (or $\omega_E$-detectable) if there exists an operator $H_\omega : \mathcal{O} \to H^1(\omega)$ such that the operator $(A - H_\omega C)$ generates a strongly continuous semi-group $(S_{H_\omega}(t))_{t \geq 0}$, which is $\omega_E$-stable.

### 3. $\Gamma_E$-DETECTABILITY AND $\Gamma_E$-OBSERVABILITY

It has been shown that a system which is exactly observable is detectable [2]. For linear systems, we recall the exactly $\Gamma$-observable [9]. Thus, regional observability definitions have been extended to regional boundary case for parabolic, hyperbolic as in [17, 19, 24] linear, semi-linear and nonlinear [26-27], and with duel concept [25]. However, in this subsection we presents some definitions which will be used to explain the notion of the regional boundary exponential detectability and observability in the state space $H^{1/2}(\Gamma)$ which is an extension from ref.s [7, 9].

**Definition 3.1:** The semi-group $(S_A(t))_{t \geq 0}$ is regionally boundary exponentially stable in $H^{1/2}(\Gamma)$ (or $\Gamma_E$-stable) if, for some positive constants $F_\Gamma$ and $\sigma_\Gamma$, then

$$\|\chi_\omega S_A(.)\|_{L(H^{1/2}(\Gamma),\ X)} \leq F_\Gamma e^{-\sigma_\Gamma t}, t \geq 0 \qquad (4)$$

In this work, we only need the relation (4) to be true on a given sub-region $\Gamma$ of the region $\partial\Omega$ in the following result.

**Remark 3.2:** If the semi-group $(S_A(t))_{t \geq 0}$ is regionally boundary exponentially stable on $H^{1/2}(\Gamma)$, then for all $x_0 \in H^{1/2}(\Gamma)$, the solution associated to the autonomous system of (1) satisfies

$$\lim_{t \to \infty} \|\gamma_0 x(.,t)\|_{H^{1/2}(\Gamma)} = \lim_{t \to \infty} \|\gamma_0 S_A x_0(.)\|_{H^{1/2}(\Gamma)} = 0 \qquad (5)$$

**Definition 3.3:** The system (1) is said to be regionally boundary exponentially stable on Γ (or $\Gamma_E$-stable), if the operator $A$ generates a semi-group which is exponentially stable on the space $H^{1/2}(\Gamma)$.

**Definition 3.4:** The system (1)-(2) is said to be regionally boundary exponentially detectable on Γ (or $\Gamma_E$-detectable) if there exists an operator

$H_\Gamma : \mathbb{R}^q \to H^{1/2}(\Gamma)$ such that $(A - H_\Gamma C)$,

generates a strongly continuous semi-group $(S_{H_\Gamma}(t))_{t \geq 0}$ which is $\Gamma_E$-stable.

However, one can deduce the following results : Thus, the notion of $\Gamma_E$-detectability is a weaker property than the exact Γ-observability as in (ref.s [3, 18]).

**Corollary 3.5:** If the systems (1)-(2) are exactly Γ-observable, then it is $\Gamma_E$-detectable. This result allows

$\exists \gamma > 0$ such that

$$\|x_\Gamma \gamma S_A x_0(.)\|_{H^{1/2}(\Gamma)} \leq v \|C S_A x_0(.)\|_{L^2(0,\infty,\mathcal{O})}, \forall x_0 \in H^{1/2}(\Gamma). \qquad (6)$$

**Proof:** We conclude the proof of this corollary from the results on observability considering $x_\Gamma K^*$ [20]. We have the following forms:

(a)   $\text{Im } F \subset \text{Im } G$
(b)   There exists $v > 0$ such that $\|F^* x^*\|_{P^*} \leq v \|G^* x^*\|_{U^*}, \forall x^* \in V^*$.

From the right hand said of this relation $\exists M, \alpha > 0$ with $v < M$ such that

$v \|G^* x^*\|_{U^*} \leq M e^{-\alpha t} \|x^*\|_{U^*}$

where $P$, $U$ and $V$ be a Banach reflexive space and $F \in L(P,V), G \in L(U,V)$.

Now, Let $P = V = H^{1/2}(\Gamma)$, $U = 0$, $F = I$ to $H^{1/2}(\Gamma)$ and $G = S^*_A(.) \chi^*_\Gamma \gamma^* C^*$ where $S_A(.)$ is a strongly continuous semi-group generates by $A$, which is $\Gamma_E$-stable then, it is $\Gamma_E$-detectable ∎.

As in El Jai and Pritchard [28], we will develop a characterization result that links the $\Gamma_E$-detectability in terms of sensors structures.

**Proposition 3.6:** Suppose that there are $q$ zone sensors $(D_i, f_i)_{1 \leq i \leq q}$. If

(1)  $q \geq r$.
(2)  $\text{Rank } G_n = r_n, \forall n, \ n = 1, .., J$

with $G = (G_n)_{ij} = (\langle \varphi_{nj}, f_i \rangle_{L^2(Di)})$ where $sup\, r_n = r < \infty$ and     $j = 1, \dots, r_n$.

Then the systems (1)-(2) are $\Gamma_E$-detectable.

**Proof:** Since rank condition is satisfied, then the systems (1)-(2) are weakly Γ-observable for finite sub-systems of (1) (see [3]) and then are exactly Γ-observable. Thus, from previous corollary 3.5 we have that the systems (1)-(3) are $\Gamma_E$-detectable.

## 4. REGIONAL EXPONENTIAL FULL ORDER OBSERVABILITY

In this section, we present an approach which allows construction an reconstruct regional exponential full order estimator ($\Gamma_{EFO}$- estimator) of $\hat{T}x(\xi,t)$. This method avoids the calculation of the inverse operators, and the consideration of the initial state [16, 19]. It enables to estimate the current state in Γ without needing the effect of the initial state of the original system.

### 4.1 Reconstruction of $\Gamma_{EFO}$-Estimator

We consider the system and the output specified by the following form:

$$\begin{cases} \frac{\partial x}{\partial t}(\xi,t) = Ax(\xi,t) + Bu(t) & Q \\ x(\xi,t) = 0 & \Theta \\ x(\xi,0) = x_0(\xi) & \Omega \\ y(.,t) = Cx(.,t) & Q \end{cases} \qquad (7)$$

Let $\Gamma \subset \partial\Omega$ be a given subdomain (region) of $\partial\Omega$ and assume that for $T \in L\left(H^{1/2}(\Gamma)\right)$, and $\hat{T} = \chi_\Gamma T$ there exists a system with state $z(.,t)$ such that

$$z(\xi,t) = \hat{T}\, x(\xi,t). \qquad (8)$$

Thus, if we can build a system which is an exponential estimator for $z(.,t)$, then it will be an exponential estimator for $\widehat{T\gamma_0}\, x(\xi,t)$, that is to say an exponential estimator to the restriction of $Tx(\xi,t)$ to the region $\Gamma$. The equations (2)-(8) give

$$\begin{bmatrix} y \\ z \end{bmatrix} = \begin{bmatrix} C \\ \hat{T} \end{bmatrix} x. \qquad (9)$$

If we assume that there exists two linear bounded operators $R$ and $S$, where $R: \mathbb{R} \to H^{1/2}(\Gamma)$, and $S: H^{1/2}(\Gamma) \to H^{1/2}(\Gamma)$, such that $RC + S\hat{T} = I$, then by deriving $z(\xi,t)$ we have

$$\frac{\partial z}{\partial t}(\xi,t) = \hat{T}\frac{\partial x}{\partial t}(\xi,t) = \hat{T}Ax(\xi,t) + \hat{T}Bu(t)$$

$$= \hat{T}ASz(\xi,t) + \hat{T}ARy(.,t) + \hat{T}Bu(t).$$

Consider now the system ($\Gamma_{EFO}$- estimator for $x$)

$$\begin{cases} \frac{\partial \hat{z}}{\partial t}(\xi,t) = F_\Gamma \hat{z}(\xi,t) + G_\Gamma u(t) + H_\Gamma y(.,t) & Q \\ \hat{z}(\xi,t) = 0 & \Theta \\ \hat{z}(\xi,0) = \hat{z}_0(\xi) & \Omega \end{cases} \qquad (10)$$

Where $F_\Gamma$ generates a strongly continuous semi-group $(S_{F_\Gamma}(t))_{t\geq 0}$, which is regionally exponential stable on $X = H^{1/2}(\Gamma)$, i.e., $\exists\, M_{F_\Gamma},\ \alpha_{F_\Gamma} > 0$, such that

$$\left\|\chi_\Gamma \gamma_0 S_{F_\Gamma}(.)\right\|_{L\left(H^{1/2}(\Gamma), H^{1/2}(\Gamma)\right)} \leq M_{F_\Gamma} e^{-\alpha_{F_\Gamma} t},\ \forall t \geq 0. \qquad (11)$$

and $G_\Gamma \in L(\mathbb{R}^p, H^{1/2}(\Gamma))$ and $H_\Gamma \in L(\mathbb{R}^q, H^{1/2}(\Gamma))$. The solution of (10) is given by

$$\hat{z}(.,t) = S_{F_\Gamma}(t)\hat{z}_0(.) + \int_0^t S_{F_\Gamma}(t-\tau)[G_\Gamma u(\tau) + H_\Gamma y(.,\tau)]d\tau \qquad (12)$$

### 4.2 $\Gamma_{EFO}$-Observability

In this case, we can consider $\hat{T} = I$, and $X = Z$, so the operator equation $\hat{T}A - F_\Gamma \hat{T} = H_\Gamma C$ of the $\Gamma_{EFO}$-observable becomes to $F_\Gamma = A - H_\Gamma C$, where $A$ and $C$ are known. Thus, the operator $H_\Gamma$ must be determined such that the operator $F_\Gamma$ is $\Gamma_{EFO}$- stable. For the system (7), consider the dynamical system

$$\begin{cases} \frac{\partial \hat{z}}{\partial t}(\xi,t) = A\hat{z}(\xi,t) + Bu(t) + H_\Gamma(y(.,t) - C\hat{z}(\xi,t) & Q \\ \hat{z}(\eta,t) = 0 & \Theta \\ \hat{z}(\xi,0) = 0 & \Omega \end{cases} \qquad (13)$$

Thus, a sufficient condition for existence of $\Gamma_{EFO}$- estimator is formulated in the following proposition.

**Proposition 4.1:** Suppose that the systems (1)-(2) are $\Gamma_{EFO}$-detectable, and then the dynamical system (13) achieve the $\Gamma_{EFO}$-observability for the systems (1)-(2), *i.e.*,

$$\lim_{t\to\infty}\|x(\xi,t)-\hat{z}(\xi,t)\|_{H^{1/2}(\Gamma)} = 0.$$

**Proof :** By the same way with minor modifications as in ref. [6] we can prove the proposition 4.1 in different case of sensors (zone, pointwise) internal or boundary.

### 4.3 Crossing Method From Internal To Boundary Case

The regional boundary exponential observability in $\Gamma$ (or $\Gamma_E$ -observability) may be seen as internal regional exponential observability in $\omega_r \subset \Omega$ (or $\omega_{r_{EFO}}$-observability) if we consider the following:

Let $\Re$ be the continuous linear extension operator [23]

$$\Re : H^{1/2}(\Gamma) \longrightarrow H^1(\Omega) \text{ such that}$$

$$\chi_\Gamma \gamma_0 \Re h(\xi,t) = h(\xi,t), \quad \forall h \in H^{1/2}(\Gamma) \tag{14}$$

Let $r > 0$ is an arbitrary and sufficiently small real and let the sets

$$E = \bigcup_{x\in\Gamma} B(z,r) \text{ and } \omega_r = E \cap \Omega \tag{15}$$

where $B(x,r)$ is the ball of radius $r$ centered in $x(\xi,t)$ and $\Gamma$ is a part of $\bar{\omega}_r$ (see Figure 3).

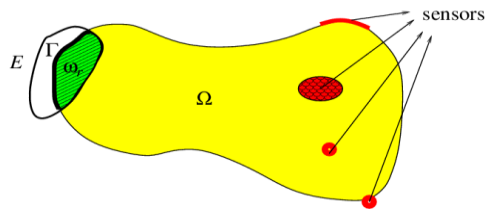

**Fig. 3:** The Domain $\Omega$, Sub-Domain $\omega_r$ and the Region $\Gamma$.

**Definition 4.2:** If The system (1) is $\omega_{r_E}$-stable, then the solution of autonomous system associated with (1)-(2), converges exponentially to zero when $t$ tends to $\infty$.

**Definition 4.3:** The systems (1)-(2) are said to be regionally exponential detectable on $\omega_r$ (or $\omega_{r_E}$-detectable) if there exists an operator $H_{\omega_r}: \mathcal{O} \to H^1(\omega_r)$ such that the operator $(A - H_{\omega_r}C)$ generates a strongly continuous semi-group $(S_{H_{\omega_r}}(t))_{t\geq 0}$, which is -stable.

Now, the method of crossing from internal $\omega_{r_E}$-detectability into $\Gamma_E$-detectability [29] will be given in a proposition below.

**Proposition 4.4:** If the systems (1)-(2) are $\bar{\omega}_{r_E}$-detectable, then it is $\Gamma_E$-detectable.

**Proof:** Let $x(\xi,t) \in H^{1/2}(\Gamma)$ and $\bar{x}(\xi,t)$ be an extension to $H^{1/2}(\partial\Omega)$. By using equation (14) and trace theorem, there exist $\Re\bar{x}(\xi,t) \in H^1(\Omega)$ [23] with a bounded support such that

$$\gamma_0(\Re\bar{x}(\xi,t)) = \bar{x}(\xi,t)$$

Since the systems (1)-(2) are $\bar{\omega}_{r_E}$-detectable, then it is $\omega_r$-detectable [19, 29]. Thus, there exists an operator $\chi_{\omega_r}K^*: \mathcal{O} \to H^1(\omega_r)$ given by

$$H_{\omega_r} x(.,t) = \chi_{\omega_r} K^* y(\xi, t)$$

such that the operator $(A - H_{\omega_r} C)$ generates a strongly continuous semi-group $(S_{\omega_r}(t))_{t \geq 0}$ which is $\omega_r$-stable. For every $\in \mathcal{O}$, we then obtain

$$\chi_{\omega_r} K^* y(\xi, t) = \chi_{\omega_r} \Re \bar{x}(\xi, t)$$

and hence

$$\chi_\Gamma (\gamma_0 \chi_{\omega_r} K^* y)(.,t) = x(\xi, t)$$

Consequently there exists an operator

$$H_\Gamma = \chi_\Gamma (\gamma_0 \chi_{\omega_r} K^* y) : \mathcal{O} \to H^{1/2}(\Gamma)$$

Such that $(A - H_\Gamma C)$ generates a semi-group $(S_{H_\Gamma}(t)) t \geq 0$ which is $\Gamma_E$-stable. Finally, the systems (1)-(2) are $\Gamma_E$-detectable. ■

**Definition 4.5:** The system (4) augmented with the output function (2) is said to be exactly observable on $\omega_r$ (or exactly $\omega_r$-observable), if

$$Im (\chi_{\omega_r} K^*) = H^1(\omega_r)$$

The following proposition shows that $\bar{\omega}_{r_{EFO}}$-observability lead $\Gamma_{EFO}$–observability.

**Proposition 4.6:** If the systems (1)-(2) are $\bar{\omega}_{r_{EFO}}$-observability, then it is $\Gamma_{EFO}$-observability.

**Proof**: By using the same hypotheses in proposition 4.4 such that

$$x(\xi, t) \in H^{1/2}(\Gamma) \text{ and } \bar{x}(\xi, t)$$

be an extension to $H^{1/2}(\partial \Omega)$. Thus, from equation (14) and trace theorem, there exist $\Re \bar{x}(\xi, t) \in H^1(\Omega)$ with a bounded support such that

$$\gamma_0 (\Re \bar{x}(\xi, t)) = \bar{x}(\xi, t)$$

Since the systems (1)-(2) are regionally exponentially full order observable on $\bar{\omega}_r$ (or $\bar{\omega}_{r_{EFO}}$-observable), so we can deduce that:

1- The systems (1)-(2) are $\omega_{r_{EFO}}$-observable, thus there exists a dynamical system with $x(\xi, t) \in X$ such that:

$$\chi_\omega T x(\xi, t) = \chi_\omega \Re \bar{x}(\xi, t)$$

Then we have

$$\chi_\Gamma (\gamma_0 \chi_\omega{}^* \chi_\omega T x)(\xi, t) = x(\xi, t) \tag{16}$$

2 - The equations (2) and (16) allow:

$$\begin{bmatrix} y \\ x \end{bmatrix}(\xi, t) = \begin{bmatrix} C \\ \chi_\Gamma(\gamma_0 \chi_\omega{}^* \chi_\omega T) \end{bmatrix} x(\xi, t),$$

and there exists two linear bounded operator $\bar{R}$ and $C$ satisfy the relation

$$\bar{R} C + \chi_\Gamma (\gamma_0 \chi_\omega{}^* \chi_\omega T) = I_\Gamma,$$

3- There exist an operator $F_{\bar{\omega}_r}$ is $\bar{\omega}_{r_{EFO}}$-observable, then it is $\Gamma_E$-stable (see [29]). Finally the systems (1)-(2) are $\Gamma_{EFO}$-observable. ∎

## 5. REGIONAL BOUNDARY EXPONENTIAL REDUCED OBSERVABILITY

In this section we need some of additional assumptions, which we explain in chapter one section 1.4 for the systems state (1)-(2).

### 5.1 Regionally Boundary Reduced System

Let us consider $X = X_1 \oplus X_2$ where $X_1$ and $X_2$ are subspace of $X$. Under the hypothesis in ref. [22], we have the dynamical system given by

$$\begin{cases} \frac{\partial x_2}{\partial t}(\xi,t) = A_{21}x_1(\xi,t) + A_{22}x_2(\xi,t) + B_2 u(t) & Q \\ x_2(\eta,t) = 0 & \Theta \\ x_2(\xi,0) = x_{2_0}(\xi) & \Omega \end{cases} \quad (17)$$

augmented with output function

$$y(.,t) = C x_1(\xi,t) \quad (18)$$

The problem consists in constructing a regional exponential estimator that enables one to estimate the unknown part $x_2(\xi,t)$ equivalent; now to define the dynamical system for (18). Thus, equations (17)-(18) allow the following system:

system:

$$\begin{cases} \frac{\partial a}{\partial t}(\xi,t) = A_{22}a(\xi,t) + [B_2 u(t) + A_{12} y(.,t)] & Q \\ a(\eta,t) = 0 & \Theta \\ a(\xi,0) = a_0(\xi) & \Omega \end{cases} \quad (19)$$

with the output function

$$\tilde{y}(.,t) = A_{12}a(.,t) \quad (20)$$

where the state $a$ in system (19) plays the role of the state $x_2$ in system (17).

### 5.2 $\Gamma_{ER}$-Observability And $\Gamma_{ER}$-Detectability

As in ref. [30] we can extend these result to the case of regional reduced orderd system from regional observability and $\Gamma_E$- detectability. In this case, the equation (7) it can be given by defining the following operator

$$\mathcal{K}: x_2 \rightarrow \mathcal{K}x_2 = A_{12} S_{A_{22}}(t) x_2 \in \mathcal{O}, \text{ then } y(.,t) = \mathcal{K}x_{2_0}(.), \text{ with}$$

the adjoint $\mathcal{K}^*: \mathcal{O} \rightarrow x_2$ such that

$$\mathcal{K}^* y^*(.,t) = \int_0^t S^*_{A_{22}}(s) A^*_{12} y^*(.,s) ds.$$

Let $\Gamma \in \partial\Omega$ and $\chi_\Gamma: H^{1/2}(\Gamma) \rightarrow H^{1/2}(\Gamma) = X_2, x_2 \rightarrow \chi_\Gamma x_2 = x_2|_\Gamma$

where $x_2|_\Gamma$ is the restriction of the state $x_2$ to $\Gamma$.

**Definition 5.1:** The systems (19)-(20) are called exactly regionally boundary exponential reduced-observable (or exactly $\Gamma_{ER}$-observable) if

$$\text{Im}\chi_\Gamma \gamma_0 \mathcal{K}^* = H^{1/2}(\Gamma) = X_2$$

**Definition 5.2:** The systems (19)-(20) are called weakly regionally boundary exponential reduced-observable (or weakly $\Gamma_{ER}$-observable) if

$$\overline{\text{Im}\chi_\Gamma \gamma_0 \mathcal{K}^*} = H^{1/2}(\Gamma) = X_2$$

This equation $\overline{\text{Im}\chi_\Gamma \gamma_0 \mathcal{K}^*}$ is equivalent to $\ker \mathcal{K} \gamma_0^* \chi_\Gamma^* = \{0\}$.

**Definition 5.3:** The suite of sensors (zone or pointwise) $(D_i, f_i)_{1 \leq i \leq q}$ are called regional boundary exponential reduced strategic sensors (or $\Gamma_{ER}$-strategic sensors if the systems (19)-(20) are weakly $\Gamma_{ER}$-observable.

**Remark 5.4:** We know the semi-group $(S_{A_{22}}(t))_{t>0}$ on Hilbert space $H^{1/2}(\Gamma)$ is said to be $\Gamma_{ER}$-stable [29], if there exists $M_{A_{22}}, \alpha_{A_{22}} > 0$ such that

$$\left\| S_{A_{22}}(t) \right\|_{H^{1/2}(\Gamma)} \leq M_{A_{22}} e^{-\alpha_{A_{22}}(t)}, \ t \geq 0 \quad (21)$$

**Remark 5.5:** The relation (21) is true on a given subdomain $\Gamma \subset \partial\Omega$, *i.e.*

$$\left\| \chi_\Gamma S_{A_{22}}(t) \right\|_{L(H^{1/2}(\Gamma), X))} \leq M_{A_{22}} e^{-\alpha_{A_{22}}(t)}, \ t \geq 0 \quad (22)$$

and then

$$\lim_{t \to \infty} \| x_2(.,t) \|_{H^{1/2}(\Gamma)} = 0$$

Now, we refer to this as regional boundary exponential reduced stability (or $\Gamma_{ER}$- stability).

**Definition 5.6:** The system (19) is said to be regional boundary exponential reduced stability (or $\Gamma_{ER}$- stable) if the operator $A_{22}$ generates a semi-group which is $\Gamma_{ER}$- stable.

**Definition 5.7:** The systems (19)-(20) are said to be regional boundary exponential reduced detectability (or $\Gamma_{ER}$-detectable) if there exists an operator $\mathcal{H}_\Gamma: \mathbb{R}^q \to H^{1/2}(\Gamma)$ such that $(A_{22} - \mathcal{H}_\Gamma A_{12})$ generates a strongly continuous semi-group $(S_{A_{22}}(t))_{t \geq 0}$, which is $\Gamma_{ER}$- stable.

From proposition 4.4, we have the dynamical system for (19)-(20) may be given by

$$\begin{cases} \frac{\partial \hat{z}}{\partial t}(\xi, t) = A_{22}\hat{z}(\xi, t) + [B_2 u(t) + A_{21} y(.,t)] + \\ \qquad \mathcal{H}_\Gamma[\tilde{y}(.,t) - A_{12}\hat{z}(\xi, t)] & Q \\ \hat{z}(\eta, t) = 0 & \Theta \\ \hat{z}(\xi, 0) = \hat{z}_0(\xi) & \Omega \end{cases} \quad (23)$$

Where $(A_{22} - \mathcal{H}_\Gamma A_{12})$ generates a strongly continuous semi-group $(S_{A_{22}}(t))_{t \geq 0}$ which is $\Gamma_{ER}$- stable on the Hilbert space $X_2 \subset X = H^{1/2}(\Gamma)$,

$$(B_2 - \mathcal{H}_\Gamma B_1) \in L(\mathbb{R}^p, X_2)$$

and

$$(A_{22}\mathcal{H}_\Gamma - \mathcal{H}_\Gamma A_{12}\mathcal{H}_\Gamma - \mathcal{H}_\Gamma A_{11} + A_{21}) \in L(\mathbb{R}^p, X_2) \text{ as in [20]}.$$

The importance of reduced $\Gamma_{ER}$- detectability is possible to define a reduced $\Gamma_{ER}$- estimator for system state may be given by the following important result.

**Theorem 5.8:** If there are $q$ sensors $(D_i, f_i)_{1 \leq i \leq q}$ and the spectrum of $A_{22}$ contains $J$ eigenvalues with non-negative real parts. The systems (19)-(20) are $\Gamma_{ER}$- detectable if and only if

1. $q \geq m_2$

2. Rank $G_{2_i} = m_{2_i}, \forall i, i = 1, \dots, J$ with

$$G_2 = G_{2_i} = \begin{cases} \langle \varphi_j(.), f_i(.) \rangle_{L^2(D_i)}, \\ \varphi_j(b_i), \end{cases}$$

where $\sup m_{2_i} = m_2 < \infty$ and $j = 1, \dots, \infty$.

**Proof :** The prove is developed to the case of zone sensors in the following stapes:

1) The system (19) can be decomposed by the projections $\mathcal{P}$ and $I - \mathcal{P}$, on two parts, unstable and stable under the assumptions of proposition 4.6, where $\mathcal{P}$ and $(I - \mathcal{P})$ are played the role of projection as $E_1, E_2$ [22]. The state vector may be given by

$$x_2(\xi, t) = [x_{2_1}(\xi, t) x_{2_2}(\xi, t)]^{tr}$$

where $x_{2_1}(\xi, t)$ is the state component of the unstable part of system (19), may be written in the form

$$\begin{cases} \frac{\partial x_{2_1}}{\partial t}(\xi, t) = A_{22_1} x_{2_1}(\xi, t) + \mathcal{P}[A_{21_1} x_{1_1}(\xi, t) + B_2 u(t)] & Q \\ x_{2_1}(\eta, t) = 0 & \Theta \\ x_{2_1}(\xi, 0) = x_{2_{1_0}}(\xi) & \Omega \end{cases} \quad (24)$$

and $x_{2_2}(\xi, t)$ is the component state of the stable part of system(19), given by

$$\begin{cases} \frac{\partial x_{2_2}}{\partial t}(\xi, t) = A_{22_2} x_{2_2}(\xi, t) + \mathcal{P}[A_{21_2} x_{1_2}(\xi, t) + B_2 u(t)] & Q \\ x_{2_2}(\eta, t) = 0 & \Theta \\ x_{2_2}(\xi, 0) = x_{2_{2_0}}(\xi) & \Omega \end{cases} \quad (25)$$

The operator $A_{22_1}$ is represented by a matrix of order $(\sum_{i=1}^{J} m_{2_i}, \sum_{i=1}^{J} m_{2_i})$ given by

$$A_{22_1} = diag[\lambda_{2_1}, \dots, \lambda_{2_1}, \dots, \lambda_{2_J}, \dots, \lambda_{2_J}] \text{ and } \mathcal{P} B_2 = [G_{2_1}^{tr}, G_{2_2}^{tr}, \dots, G_{2_J}^{tr}]$$

From condition (2) of this theorem, then the suite of sensors $(D_i, f_i)_{1 \leq i \leq q}$ is $\Gamma_{ER}$- strategic for the unstable part of the system (19), the subsystem (24) is weakly regionally boundary reduced-observable in $\Gamma$ (or weakly $\Gamma_{ER}$- observable ) and since it is finite dimensional, then it is exactly regionally boundary reduced-observable in $\Gamma$ (or exactly $\Gamma_{ER}$-observable).

Therefore it is $\Gamma_{ER}$- detectable, and hence there exists an operator $\mathcal{H}_\Gamma^1$ such that $(A_{22_1} - \mathcal{H}_\Gamma^1 A_{12_1})$ which satisfies the following:

$$\exists M_\Gamma^1, \alpha_\Gamma^1 > 0 \text{ such that } \left\| e^{(A_{22_1} - \mathcal{H}_\Gamma^1 A_{12_1}) t} \right\|_{H^{1/2}(\Gamma)} \leq M_\Gamma^1 e^{-\alpha_\Gamma^1(t)}$$

and we have

$$\left\| x_{2_1}(\xi, t) \right\|_{H^{1/2}(\Gamma)} \leq M_\Gamma^1 e^{-\alpha_\Gamma^1(t)} \left\| \mathcal{P} x_{2_0}(.) \right\|_{H^{1/2}(\Gamma)}$$

Since the semi-group generated by the operator $A_{22_2}$ is $\Gamma_{ER}$-stable, $\exists M_\Gamma^2, \alpha_\Gamma^2 > 0$ such that

$$\left\| x_{2_2}(\xi, t) \right\|_{H^{1/2}(\Gamma)} \leq M_\Gamma^2 e^{-\alpha_\Gamma^2(t)} \left\| (I - \mathcal{P}) x_{2_0}(.) \right\|_{H^{1/2}(\Gamma)}$$

$$+ \int_0^t M_\Gamma^2 e^{-\alpha_\Gamma^2(t-\tau\tau)} \|(I-\mathcal{P})x_{2_0}(.)\|_{H^{1/2}(\Gamma)} \|u(t)\| d\tau$$

and there fore $x_2(\xi,t) \to 0$ when $t \to \infty$. Thus, the systems (19)-(20) are $\Gamma_{ER}$-detectable.

2) If the systems (19)-(20) are $\Gamma_{ER}$-detectable, then $\exists \mathcal{H}_\Gamma \in L(L^2(0,\infty,\mathbb{R}^q), H^{1/2}(\Gamma))$ such that $(A_{22} - \mathcal{H}_\Gamma A_{12})$ generates an $\Gamma_{ER}$-stable, strongly continuous semi-group $(S_{A_{22}}(t))_{t \geq 0}$ on the space $H^{1/2}(\Gamma)$ which satisfies the following

$$\exists M_\Gamma, \alpha_\Gamma > 0 \text{ such that } \|\chi_\Gamma S_{A_{22}}(t)\|_{H^{1/2}(\Gamma)} \leq M_\Gamma e^{-\alpha_\Gamma(t)}$$

Thus the unstable subsystem (24) is $\Gamma_{ER}$-detectable. Since this subsystem is of finite dimensional, then it is exactly $\Gamma_{ER}$-observable. Therefore (24) is weakly $\Gamma_{ER}$-observable and hence it is reduced $\Gamma_{ER}$-strategic, i.e.

$$[\mathcal{K}x_\Gamma^* x_2^*(.,t) = 0 \Longrightarrow x_2^*(.,t) = 0]. \text{ For } x_2^*(.,t) \in H^{1/2}(\Gamma)$$

We have

$$[\mathcal{K}x_\Gamma^* x_2^*(.,t) = (\sum_{j=1}^J e^{\lambda_j t} \langle \varphi_j(.), x_2^*(.,t) \rangle_{H^{1/2}(\Gamma)} \langle \varphi_j(.), f_i(.) \rangle_{H^{1/2}(\Gamma)})_{1 \leq i \leq q}$$

If the unstable system (24) is not $\Gamma_{ER}$-strategic, $\exists x_2^*(.,t) \in H^{1/2}(\Gamma)$ such that $\mathcal{K}x_\Gamma^* x_2^*(.,t) = 0$, this leads to

$$\sum_{j=1}^J \langle \varphi_j(.), x_2^*(.,t) \rangle_{H^{1/2}(\Gamma)} \langle \varphi_j(.), f_i(.) \rangle_{H^{1/2}(\Gamma)} = 0$$

the state vectors $x_{2_i}$ may be given

$$x_{2_i}(.,t) = [\langle \varphi_j(.), x_2^*(.,t) \rangle_{H^{1/2}(\Gamma)} \langle \varphi_j(.), x_2^*(.,t) \rangle_{H^{1/2}(\Gamma)}]^{tr} \neq 0$$

We then obtain $G_{2_i} x_{2_i} = 0, \forall i, i = 1, \dots, J$ and there fore $\text{Rank } G_{2_i} \neq m_{2_i}$.

Here, we construct the $\Gamma_{ER}$- estimator for parabolic distributed parameter system (1), we need to present the following remarks

**Remark 5.9:** Now, choose the following decomposition:

$$\hat{z} = \begin{bmatrix} \hat{z}_1 \\ \hat{z}_2 \end{bmatrix} = \begin{bmatrix} y \\ \varphi + \mathcal{H}_\Gamma y \end{bmatrix}$$

which estimates exponentially the state vector

$$x = \begin{bmatrix} x_1 \\ x_2 \end{bmatrix}$$

then, the dynamical system (23) is given by the following system:

$$\begin{cases} \frac{\partial \varphi}{\partial t}(\xi,t) = (A_{22} - \mathcal{H}_\Gamma A_{12})\varphi(\xi,t) + \\ \qquad [A_{22}\mathcal{H}_\Gamma - \mathcal{H}_\Gamma A_{12}\mathcal{H}_\Gamma - \mathcal{H}_\Gamma A_{11} + A_{21}] \\ \qquad y(\xi,t) + [B_2 - \mathcal{H}_\Gamma B_1]u(t) & Q \\ \varphi(\eta,t) = 0 & \Theta \\ \varphi(\xi,0) = \varphi_0(\xi) & \Omega \end{cases} \quad (26)$$

which defines an $\Gamma_{ER}$- estimator for $T_\Gamma x_2(\xi,t)$ if

1. $\lim_{t \to \infty} \|\varphi(\xi,t) - T_\Gamma x_2(\xi,t)\|_{H^{1/2}(\Gamma)} = 0$
2. $T_\Gamma: D(A_{22}) \to D(A_{22} - \mathcal{H}_\Gamma A_{12})$ where $T_\Gamma = \chi_\Gamma \gamma_0 T$ and $\varphi(\xi,t)$ is the solution of system (26).

**Remark 5.10:** The dynamical system (26) estimates the regional boundary exponential reduced state of the system (1) if the following conditions satisfies:

1. $\exists L \in L(0, H^{1/2}(\Gamma))$ and $M \in L(H^{1/2}(\Gamma))$ such that:
   $LA_{12} + MT_\Gamma = I_\Gamma$
2. $T_\Gamma A_{22} - (A_{22} - \mathcal{H}_\Gamma A_{12})T_\Gamma = \mathcal{H}_\Gamma A_{12}$ and $(B_2 - \mathcal{H}_\Gamma B_1) = T_\Gamma B_2$
3. The system (26) defines an $\Gamma_{ER}$- estimator for the system (1).
4. If $X = X_2$ and $T_\Gamma = I_\Gamma$ then, in the above case, we have
   $A_{22} - (A_{22} - \mathcal{H}_\Gamma A_{12}) = \mathcal{H}_\Gamma A_{12}$

**Remark 5.11:** The system (1) is $\Gamma_{ER}$-observable if there exists an $\Gamma_{ER}$-estimators (26) which estimates the regional boundary exponential reduced state the system. Now, we present the sufficient condition of the regional boundary exponential reduced observability notion as in the following main result.

**Theorem 5.12:** If the systems (19)-(20) are $\Gamma_{ER}$- detectable, then it is $\Gamma_{ER}$- observable by the dynamical system (25), that means

$$\lim_{t\to\infty}\left\|(\varphi(\xi,t) + \mathcal{H}_\Gamma y(\xi,t)) - x_2(\xi,t)\right\|_{H^{1/2}(\Gamma)} = 0,$$

**Proof:** The solution of the dynamical system (23) is given by

$$\hat{z}(\xi,t) = S_{\mathcal{H}_\Gamma}(t)\hat{z}_0(\xi) + \int_0^t S_{\mathcal{H}_\Gamma}(t-\tau)[B_2 u(\tau)$$
$$+ A_{21}y(\xi,\tau) + \mathcal{H}_\Gamma \tilde{y}(\xi,\tau)]d\tau \qquad (27)$$

From the equation (20), we have

$$\tilde{y}(\xi,t) = A_{12}a(.,t) = \frac{\partial x_1}{\partial t}(\xi,t) - A_{11}x_1(\xi,t) - B_1 u(t) \qquad (28)$$

By using (27) and (28), we obtain

$$\hat{z}(\xi,t) = S_{\mathcal{H}_\Gamma}(t)\hat{z}_0(\xi) + \int_0^t S_{\mathcal{H}_\Gamma}(t-\tau)\mathcal{H}_\Gamma \frac{\partial x_1}{\partial t}(\xi,t)d\tau +$$
$$\int_0^t S_{\mathcal{H}_\Gamma}(t-\tau)[B_2 u(\tau) + A_{21}y(\xi,\tau)]$$
$$-\mathcal{H}_\Gamma A_{11}x_1(.,\tau) - \mathcal{H}_\Gamma B_1 u(\tau)d\tau. \qquad (29)$$

and we can get

$$\int_0^t S_{\mathcal{H}_\Gamma}(t-\tau)\mathcal{H}_\Gamma \frac{\partial x_1}{\partial t}(\xi,t)d\tau = \mathcal{H}_\Gamma x_1(.,t) - S_{\mathcal{H}_\Gamma}(t)\mathcal{H}_\Gamma x_{1_0}(.)$$
$$+ (A_{22} - \mathcal{H}_\Gamma A_{12})\int_0^t S_{\mathcal{H}_\Gamma}(t-\tau)\mathcal{H}_\Gamma x_1(.,\tau)d\tau \qquad (30)$$

Using Bochner integrability properties and closeness of $(A_{22} - \mathcal{H}_\Gamma A_{12})$, the equation (30) becomes

$$\int_0^t S_{\mathcal{H}_\Gamma}(t-\tau)\mathcal{H}_\Gamma \frac{\partial x_1}{\partial t}(\xi,t)d\tau = \mathcal{H}_\Gamma x_1(.,t) - S_{\mathcal{H}_\Gamma}(t)\mathcal{H}_\Gamma x_{1_0}(.)$$
$$+ \int_0^t S_{\mathcal{H}_\Gamma}(t-\tau)(A_{22} - \mathcal{H}_\Gamma A_{12})\mathcal{H}_\Gamma x_1(\xi,\tau)d\tau \qquad (31)$$

Substituting (31) into (29), we have

$$\hat{z}(.,t) = S_{\mathcal{H}_\Gamma}(t)\hat{z}_0(.) - S_{\mathcal{H}_\Gamma}(t)\mathcal{H}_\Gamma x_{1_0}(.) + \mathcal{H}_\Gamma x_1(.,t)$$
$$+ \int_0^t S_{\mathcal{H}_\Gamma}(t-\tau)[A_{22}\mathcal{H}_\Gamma - \mathcal{H}_\Gamma A_{12}\mathcal{H}_\Gamma - \mathcal{H}_\Gamma A_{11} + A_{21}]$$

$$x_1(.,t)d\tau + \int_0^t S_{\mathcal{H}_\Gamma}(t-\tau)[B_2 - \mathcal{H}_\Gamma B_1]u(\tau)d\tau. \tag{32}$$

Setting $\varphi(.,t) = \hat{z}(.,t) - \mathcal{H}_\Gamma y(.,t)$, with $\varphi_0(.,0) = \hat{z}_0(.) - \mathcal{H}_\Gamma x_{1_0}(.)$, where $y_0(.) = x_{1_0}(.)$. Now, assume that $(A_{22}\mathcal{H}_\Gamma - \mathcal{H}_\Gamma A_{12}\mathcal{H}_\Gamma - \mathcal{H}_\Gamma A_{11} + A_{21})$ and $(B_2 - \mathcal{H}_\Gamma B_1)$ are bounded operators, the equation (32) can be differentiated to yield the following system

$$\begin{cases} \frac{\partial \varphi}{\partial t}(\xi,t) = (A_{22} - \mathcal{H}_\Gamma A_{12})\varphi(\xi,t) + (A_{22}\mathcal{H}_\Gamma - \mathcal{H}_\Gamma A_{12}\mathcal{H}_\Gamma \\ \quad -\mathcal{H}_\Gamma A_{11} + A_{11})y(.,t) + (B_2 - \mathcal{H}_\Gamma B_1)u(t) & Q \\ \varphi(\eta,t) = 0 & \Theta \\ \varphi(\xi,0) = \varphi_0(\xi) & \Omega \end{cases}$$

and therefore

$$\frac{\partial z}{\partial t}(\xi,t) - \frac{\partial x_2}{\partial t}(\xi,t) = (\varphi(\xi,t) + \mathcal{H}_\Gamma y(\xi,t) - x_2(\xi,t)$$

$$= (A_{22}\hat{z}(\xi,t) + B_2 u(t) + A_{21}y(.,t) + \mathcal{H}_\Gamma(\tilde{y}(\xi,t)$$

$$-A_{12}\hat{z})(\xi,t) - A_{21}x_1(\xi,t) - A_{22}x_2(\xi,t) - B_2 u(t)$$

$$= (A_{22} - \mathcal{H}_\Gamma A_{12})(\hat{z}(\xi,t) - x_2(\xi,t)) \tag{33}$$

From the relation

$$\|\chi_\Gamma S_{\mathcal{H}_\Gamma}(t)x_{2_0}(.)\|_{H^{1/2}(\Gamma)} \leq M_{\mathcal{H}_\Gamma} e^{-\alpha_{\mathcal{H}_\Gamma}(t)}$$

we obtain

$$\|\hat{z}(.,t) - x_2(.,t)\|_{H^{1/2}(\Gamma)} \leq \|\chi_\Gamma S_{\mathcal{H}_\Gamma}(t)\|_{H^{1/2}(\Gamma)}$$

$$\|\hat{z}(.,0) - x_2(.,0)\|_{H^{1/2}(\Gamma)}$$

$$\leq M_{\mathcal{H}_\Gamma} e^{-\alpha_{\mathcal{H}_\Gamma}(t)} \|\hat{z}(.,0) - x_2(.,0)\|_{H^{1/2}(\Gamma)}$$

$$\to 0 \text{ as } t \to \infty \tag{34}$$

Where the component $\hat{z}(\xi,t)$ is an exponentially estimator of $x_2$. Then, we have the system (23) is $\Gamma_{ER}$-observable for the systems (18)-(19).∎

From the previous theorem 5.12, we can deduce the following definition which characterizes another new strategic sensor:

**Definition 5.13:** A sensor is $\Gamma_{ER}$-strategic sensor if the corresponding system is $\Gamma_{ER}$-observable.

## 6. APPLICATION TO EXCHANGE SYSTEM

Consider the case of two-phase exchange systems described by the following coupled parabolic equations as in [4]

$$\begin{cases} \frac{\partial x_1}{\partial t}(\xi_1,\xi_2,t) = \alpha \frac{\partial^2 x_1}{\partial \xi^2}(\xi_1,\xi_2,t) + \beta(x_1(\xi_1,\xi_2,t) - x_2(\xi_1,\xi_2,t))Q \\ \frac{\partial x_2}{\partial t}(\xi_1,\xi_2,t) = \gamma \frac{\partial^2 x_2}{\partial \xi^2}(\xi_1,\xi_2,t) + \beta(x_1(\xi_1,\xi_2,t) - x_2(\xi_1,\xi_2,t))Q \\ x_1(\eta_1,\eta_2,t) = 0, \ x_2(\eta_1,\eta_2,t) = 0 & \Theta \\ x_1(\xi_1,\xi_2,t) = x_{1_0}(\xi_1,\xi_2), \ x_1(\xi_1,\xi_2,t) = x_{1_0}(\xi_1,\xi_2) & \Omega \end{cases} \tag{35}$$

and consider $\Omega = (0,1) \times (0,1)$ with sub-region $\Gamma = (\alpha_1, \beta_1) \times (\alpha_2, \beta_2) \subset \partial\Omega$. Suppose that it is possible to measure the states $x_1(.,t)$, by using $q$ zone sensor$(D_i, f_i)_{1 \leq i \leq q}$. The output function (2) is given by

$$y(t) = Cx_1(.,t)$$
$$= \left[\int_{D_1} x_1(\xi_1, \xi_2, t) f_1(\xi_1, \xi_2) d\xi_1 \xi_2, \ldots, \int_{D_q} x_1(\xi_1, \xi_2, t) f_q(\xi_1, \xi_2) d\xi_1 \xi_2\right]^{tr}$$

Now, the problem is to estimate exponentially $x_1(\xi_1, \xi_2, t)$. Consider now

$$\frac{\partial x}{\partial t} = \begin{bmatrix} \frac{\partial x_1}{\partial t} \\ \frac{\partial x_2}{\partial t} \end{bmatrix} = \begin{bmatrix} A_{11} & A_{12} \\ A_{21} & A_{22} \end{bmatrix} \begin{bmatrix} x_1 \\ x_2 \end{bmatrix} \qquad (36)$$

where

$$A_{11} = \alpha \frac{\partial^2 x_1}{\partial \xi^2}(\xi_1, \xi_2, t) + \beta, \ A_{22} = \gamma \frac{\partial^2 x_1}{\partial \xi^2}(\xi_1, \xi_2, t) + \beta$$

and $A_{12} = A_{21} = -\beta I$.

From theorem 5.12, we can construct regional boundary reduced estimator for system (35) if the sensors$(D_i, f_i)_{1 \leq i \leq q}$ are $\Gamma$-strategic for the unstable part of the subsystem

$$\begin{cases} \frac{\partial x_1}{\partial t}((\xi_1, \xi_2, t) = \gamma \frac{\partial^2 x_1}{\partial \xi^2}(\xi_1, \xi_2, t) + \beta(x_1(\xi_1, \xi_2, t) \\ \qquad\qquad -x_1(\xi_1, \xi_2, t) & Q \\ x_1(\eta_1, \eta_2, t) = 0 & \Theta \\ x_1(\xi_1, \xi_2, 0) = x_{1_0}(\xi_1, \xi_2) & \Omega \end{cases} \qquad (37)$$

where that $\gamma = 0.1$ and $\beta = 1$. If we choose the sensors $(D_i, f_i)_{1 \leq i \leq q}$ such that

$$y(t) = \left[\int_{D_1} x_1(\xi_1, \xi_2, t) f_1(\xi_1, \xi_2) d\xi_1 \xi_2, \ldots, \int_{D_q} x_1(\xi_1, \xi_2, t) f_q(\xi_1, \xi_2) d\xi_1 \xi_2\right]^{tr} \neq 0,$$

then, there exists $\mathcal{H}_\Gamma \in L(\mathbb{R}^q, H^{\frac{1}{2}}(\Gamma))$ such that the operator $(A_{22} - \mathcal{H}_\Gamma A_{12})$ generates a strongly continuous stable group on the space $H^{1/2}(\Gamma)$. Thus we have

$$\lim_{n \to \infty} \|w(.,t) + \mathcal{H}_\Gamma x_1(.,t) - x_2(.,t)\|_{H^{1/2}(\Gamma)} = 0,$$

where

$$\begin{cases} \frac{\partial w}{\partial t}(\xi_1, \xi_2, t) = \gamma \frac{\partial^2 w}{\partial \xi^2}(\xi_1, \xi_2, t) + \beta(1 + \mathcal{H}_\Gamma) w(\xi_1, \xi_2, t) \\ \qquad\qquad + (\gamma - \alpha \mathcal{H}_\Gamma) \frac{\partial^2 x_1}{\partial \xi^2}(\xi_1, \xi_2, t) + \beta (\mathcal{H}_\Gamma^2 - 1)(\xi_1, \xi_2, t) & Q \\ w(\eta_1, \eta_2, t) = 0 & \Theta \\ w(\xi_1, \xi_2, 0) = w_0(\xi_1, \xi_2) & \Omega \end{cases} \qquad (38)$$

In this section, we give the specific results related to some examples of sensors locations and we apply these results to different situations of the domain, which usually follow from symmetry considerations.

We consider the two-dimensional system defined on $\Omega = (0,1) \times (0,1)$ with the case of system described by the following equations:

$$\begin{cases} \frac{\partial x_2}{\partial t}(\xi_1,\xi_2,t) = \gamma \frac{\partial^2 x_2}{\partial \xi^2}(\xi_1,\xi_2,t) + \beta(x_2(\xi_1,\xi_2,t) \\ \qquad\qquad - \beta x_1(\xi_1,\xi_2,t) & Q \\ x_2(\eta_1,\eta_2,t) = 0 & \Theta \\ x_2(\xi_1,\xi_2,0) = x_{2_0}(\xi_1,\xi_2) & \Omega \end{cases} \quad (39)$$

augmented with the output function

$$y(t) = Cx_1(.,t) \quad (40)$$

Let $\Gamma = \Pi_{i=1}^2(\alpha_1,\beta_1) \times (\alpha_2,\beta_2)$, In this case the eigenfunctions and eigenvalues for the dynamic system (39) given by

$$\varphi_{ij}(\xi_1,\xi_2) = \frac{2}{\sqrt{(\beta_1-\alpha_1)(\beta_2-\alpha_2)}} \sin i\pi\left(\frac{(\xi_1-\alpha_1)}{(\beta_1-\alpha_1)}\right) \sin j\pi\left(\frac{(\xi_2-\alpha_2)}{(\beta_2-\alpha_2)}\right) \quad (41)$$

$$\lambda_{ij} = -\left(\frac{i^2}{(\beta_1-\alpha_1)^2} + \frac{j^2}{(\beta_2-\alpha_2)^2}\right)\pi^2, \quad i,j \geq 1 \quad (42)$$

we examine the two cases illustrated in (Figures 18 and 19).

### 6.1 Internal Rectangular Zone Sensor

For discussing this case, suppose the systems (39)-(40) where the sensor supports $D_i$ is the located in $\Omega$ as in (Figure 4). The output function can be written by the form

$$y(t) = \int_{D_1} x_1(\xi_1,\xi_2,t) f_1(\xi_1,\xi_2) d\xi_1 \xi_2, \quad (43)$$

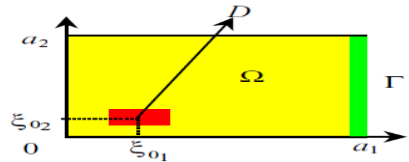

**Fig. 4:** Rectangular Domain, Region $\Gamma$ and Location $D$ with Rectangular Support Sensor.

Then, the sensor $(D_i,f_i)_{1\leq i\leq q}$ may be sufficient for $\Gamma_{ER}$-observability, and there exists $\mathcal{H}_\Gamma \in L(\mathbb{R}^q, H^{\frac{1}{2}}(\Gamma))$ such that the operator $(A_{22} - \mathcal{H}_\Gamma A_{12})$ generates a strongly continuous stable semi-group on the space $H^{1/2}(\Gamma)$. Thus we have

$$\lim_{t\to\infty} \|(w(\xi_1,\xi_2,t) + \mathcal{H}_\Gamma x_2(\xi_1,\xi_2 t)) - x_1(\xi_1,\xi_2,t)\|_{H^{1/2}(\Gamma)} = 0,$$

where

$$\begin{cases} \frac{\partial w}{\partial t}((\xi_1,\xi_2,t) = \gamma \frac{\partial^2 w}{\partial \xi^2}(\xi_1,\xi_2,t) + \beta(1+\mathcal{H}_\Gamma)w(\xi_1,\xi_2,t) \\ \qquad + (\gamma - \alpha\mathcal{H}_\Gamma)\frac{\partial^2 x_1}{\partial \xi^2}(\xi_1,\xi_2,t) + \beta\,(\mathcal{H}_\Gamma^2 - 1)(\xi_1,\xi_2,t) & Q \\ w(\eta_1,\eta_2,t) = 0 & \Theta \\ w(\xi_1,\xi_2,0) = w_0(\xi_1,\xi_2) & \Omega \end{cases} \quad (44)$$

If $D_i = \Pi_{i=1}^2[\xi_{0_i} - \iota_i, \xi_{0_i} + \iota_i]$, with $[\xi_{0_i} - \alpha_i/\xi_{0_i} + \alpha_i] \in Q$ then, we have the following result.

**Proposition 6.1:** Let $f_i$ are symmetric about line $x_{0_i} = \xi_{0_i}$ and the sensors $(D_i,f_i)_{1\leq i\leq q}$ are not strategic for the systems (39)-(40), and then these systems are not $\Gamma_{ER}$-observable by the $\Gamma_{ER}$-estimator systems (3.44). If for any $i_0 \in 1 \leq i \leq q$ such that

$$\frac{i_0(\xi_{0_i}-\alpha_i)}{\beta_1-\alpha_1}, \frac{j_0(\xi_{0_i}-\alpha_i)}{\beta_2-\alpha_2} \in Q$$

**Proof:** Suppose that $i_0 = 1$, and $[\beta_{1-}\alpha_1/\beta_{2-}\alpha_2] \in Q$, then there exists $j_0 \geq 1$ such that

$Sin\ j_0(\beta_{1-}\alpha_1/\beta_{2-}\alpha_2) = 0$. But

$$y(t) = \langle f_1, \varphi_{0_1} j_0 \rangle = \left(\frac{4}{(\beta_1-\alpha_1)(\beta_2-\alpha_2)}\right)^{1/2} \int_{\alpha_2-l_2}^{\alpha_2+l_2} \int_{\alpha_1-l_1}^{\alpha_1+l_1} f_1(\xi_1, \xi_2)$$

$$\sin\left[\frac{j_0\pi(\xi_{0_1}-\alpha_1)}{\beta_1-\alpha_1}\right] \sin\left[\frac{j_0\pi(\xi_{0_2}-\alpha_2)}{(\beta_2-\alpha_2)}\right] d\xi_1 d\xi_2 = 0$$

**6.2 Internal Pointwise Sensor**

Consider the case of pointwise sensor located inside of $\Omega$. The system (39) augmented with the following output function:

$$y(t) = \int_{\Omega\backslash\omega} x_2(\xi_1, \xi_2, t) \delta(\xi_1 - b_1, \xi_2 - b_2) d\xi_1 d\xi_2 \qquad (45)$$

Where $b = (b_1, b_2) \in \Omega$ as in (Figure 5) is the location of pointwise sensor with $(b_1 - \alpha_1)/(\beta_1 - \alpha_1)$ and $(b_2 - \alpha_2)/(\beta_2 - \alpha_2) \in Q$.

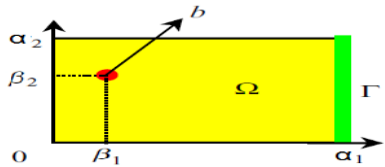

**Fig. 5:** Rectangular Domain $\Omega$, Region $\Gamma$.

Then, the sensor $(b, \delta_b)$ may be sufficient for $\Gamma_{ER}$-observability [9], and there exists $\mathcal{H}_\Gamma \in L(\mathbb{R}^q, H^{\frac{1}{2}}(\Gamma))$ such that the operator $(A_{22} - \mathcal{H}_\Gamma A_{12})$ generates a strongly continuous stable semi-group on the space $H^{1/2}(\Gamma)$. Thus we have

$$\lim_{t\to\infty} \left\|\left(w(\xi_1, \xi_2, t) + \mathcal{H}_\Gamma x_2(\xi_1, \xi_2 t)\right) - x_1(\xi_1, \xi_2, t)\right\|_{H^{1/2}(\Gamma)} = 0,$$

where

$$\begin{cases} \frac{\partial w}{\partial t}(\xi_1, \xi_2, t) = \gamma \frac{\partial^2 w}{\partial \xi^2}(\xi_1, \xi_2, t) + \beta(1 + \mathcal{H}_\Gamma)w(\xi_1, \xi_2, t) \\ \qquad + (\gamma - \alpha\mathcal{H}_\Gamma)\frac{\partial^2 x_1}{\partial \xi^2}(\xi_1, \xi_2, t) + \beta(\mathcal{H}_\Gamma^2 - 1)(\xi_1, \xi_2, t) & Q \\ w(\eta_1, \eta_2, t) = 0 & \Theta \\ w(\xi_1, \xi_2, 0) = w_0(\xi_1, \xi_2) & \Omega \end{cases} \qquad (46)$$

Then, we have the following result

**Corollary 6.2:** The systems (39)-(45) are not $\Gamma_{ER}$-observable by the $\Gamma_{ER}$-estimator (46), If for any $i_0 \in 1 \leq i \leq 2$, $j_0 \in 1 \leq i \leq q$ such that

$$\frac{i_0(b_1-\alpha_1)}{\beta_1-\alpha_1}, \frac{j_0(b_2-\alpha_2)}{\beta_2-\alpha_2} \in Q$$

**Proof:** Assume that $i_0 = 1$, and $[\beta_{1-}\alpha_1/\beta_{2-}\alpha_2] \in Q$, then there exists $j_0 \geq 1$ such that

$Sin\ j_0(b_1 - \alpha_1)/(\beta_1 - \alpha_1) = 0$. But

$$y(t) = \langle f_1, \varphi_{0_1} j_0 \rangle = \left(\frac{4}{(\beta_1-\alpha_1)(\beta_2-\alpha_2)}\right)^{1/2} \int_{\Omega/\omega} \delta_b(b_1, b_2)$$

$$\sin\left[\frac{j_0\pi(b_1-\alpha_1)}{\beta_1-\alpha_1}\right] \sin\left[\frac{j_0\pi(b_2-\alpha_2)}{(\beta_2-\alpha_2)}\right] d\xi_1 d\xi_2 = 0$$

**Remark 6.3:** These results can be extended to the following:

(1) Case of Neumann or mixed boundary conditions [29].

(2) Case of boundary (pointwise, filament and zone) sensors as in [31].

**CONCLOSION**

The concept developed in this paper is related to the regional boundary exponential reduced observability in connection with the strategic sensors characterizations. Various interesting results concerning the choice of sensors are given and illustrated in specific situations. Many questions still opened. For example, the problem of finding the optimal sensor location ensuring such an objective. The result of regional exponential reduced observability concept of hyperbolic linear or semi linear or nonlinear systems is under consideration.

**ACKNOWLEDGMENTS:**

Our thanks in advance to the editors and experts for considering this paper to publish in this estimated journal and for their efforts.